\documentclass[a4paper,12pt]{article}

\usepackage{color}

\usepackage{makeidx}

\usepackage{latexsym}
\usepackage{amscd}
\usepackage{amsmath}
\usepackage{amssymb}

\usepackage{amsthm}
\usepackage{float}
\usepackage{graphicx}
\usepackage{psfrag}

\theoremstyle{plain}
\newtheorem{theorem}{Theorem}
\newtheorem{corollary}[theorem]{Corollary}
\newtheorem{lemma}[theorem]{Lemma}

\newtheorem{proposition}[theorem]{Proposition}

\theoremstyle{definition}
\newtheorem{definition}[theorem]{Definition}
\newtheorem{remark}[theorem]{Remark}
\newtheorem{example}[theorem]{Example}

\newdimen\argwidth
\def\db[#1\db]{%
 \setbox0=\hbox{$#1$}\argwidth=\wd0
 \setbox0=\hbox{$\left[\box0\right]$}
  \advance\argwidth by -\wd0
 \left[\kern.3\argwidth\box0 \kern.3\argwidth\right]}

\newcommand{\codim}{\operatorname{codim}}

\newcommand{\bC}{\ensuremath{\mathbb{C}}}

\newcommand{\bR}{\ensuremath{\mathbb{R}}}

\newcommand{\bZ}{\ensuremath{\mathbb{Z}}}

\newcommand{\scA}{\ensuremath{\mathcal{A}}}

\newcommand{\scF}{\ensuremath{\mathcal{F}}}

\newcommand{\scL}{\ensuremath{\mathcal{L}}}

\newcommand{\scS}{\ensuremath{\mathcal{S}}}

\newcommand{\sfch}{\mathsf{ch}}
\newcommand{\sfM}{\mathsf{M}}

\newcommand{\ch}{\mathsf{ch}}
\newcommand{\bch}{\mathsf{bch}}
\newcommand{\Poin}{\operatorname{Poin}}

\newcommand{\ii}{\sqrt{-1}}
\newcommand{\sign}{\operatorname{sgn}}
\newcommand{\wt}{\widetilde}
\newcommand{\ddt}{\left.\frac{d}{dt}\right|_{t=0}}
\newcommand{\owari}{\hfill$\square$\medskip}
\title{Chamber basis of the Orlik-Solomon algebra and Aomoto complex}
\author{Masahiko Yoshinaga}
\date{\today}
\begin{document}
\maketitle

\begin{abstract}
We introduce a basis of 
the Orlik-Solomon algebra labeled by chambers, 
so called chamber basis. 
We consider structure constants of the 
Orlik-Solomon algebra 
with respect to the chamber basis and 
prove that these structure constants recover 
D. Cohen's minimal complex from the Aomoto complex.

\noindent
{\bf MSC-class}: 32S22 (Primary) 52C35, 32S50 (Secondary) 
\end{abstract}

\section{Introduction}

Let $\scA=\{H_1, \ldots, H_n\}$ be an affine hyperplane 
arrangement in the real vector space $\bR^\ell$. 
Choose for each $H\in\scA$ an affine linear form 
$\alpha_H$ with $H=\alpha_H^{-1}(0)$. 
Denote by $\ch(\scA)$ the set of all chambers and by 
$\sfM(\scA)=\bC^\ell\backslash\bigcup_{H\in\scA}H\otimes\bC$ 
the complement to the complexified hyperplanes. 

The set $\ch(\scA)$ of chambers has been known to 
carry information about topology of $\sfM(\scA)$. 
For example $|\ch(\scA)|=\sum_{i=0}^{\ell}b_i(\sfM(\scA))$ 
\cite{zas-face}, the homotopy type of $\sfM(\scA)$ 
can be obtained from the face poset \cite{sal-top}, 
and \cite{koh-hom} uses bounded chambers to 
construct a basis of local system cohomology group. 
The relation between Orlik-Solomon algebra and 
the ring $\bZ[\ch(\scA)]$ of $\bZ$-valued 
functions over the set of chambers was studied 
in \cite{var-gel}. \cite{yos-min} and \cite{sal-sett} 
considered the 
relation between structures of chambers and 
minimal CW-decomposition. 
We will pursue these topological 
interpretations of chambers in the 
context of rank one local system cohomology groups.

Let $\lambda_H\in\bC$ be complex weights. A rank one 
local system $\scL$ on $\sfM(\scA)$ is defined with 
monodromy $\exp(2\pi\ii\lambda_H)$ around the 
hyperplane $H$. Let 
$\omega_H=\frac{1}{2\pi\ii}\frac{d\alpha_H}{\alpha_H}$, 
$A^\bullet=H^\bullet(\sfM(\scA), \bC)$ and 
$$
\omega_\lambda=\sum_{H\in\scA}\lambda_H\omega_H. 
$$
Under some genericity conditions on the weights 
$\lambda_H$, Esnault-Schechtman-Viehweg \cite{esv} proved 
that the Aomoto complex $(A^\bullet, 2\pi\ii\omega_\lambda\wedge)$ 
is quasi-isomorphic to the de Rham complex with 
coefficients in $\scL$. In particular, if the weights 
$\lambda_H$ are 
sufficiently small, e.g. $|\lambda_H|<\frac{1}{2(n+1)}\ 
(\forall H\in\scA)$, 
then 
$$
H^p(\sfM(\scA), \scL)\cong 
H^p(A^\bullet, 2\pi\ii\omega_\lambda\wedge). 
$$
Thus the local system cohomology group 
$H^\bullet(\sfM(\scA), \scL)$ can be calculated from 
the Aomoto complex if $\scL$ is 
close to the trivial one (tangent-cone theorem). 
However, in general, 
these two cohomology groups have different dimensions 
\cite{coh-suc, suc-tra}. The Aomoto complex does not 
compute $H^\bullet(\sfM(\scA), \scL)$ at least directly. 
This suggests the problem whether the 
Aomoto complex 
$(A^\bullet, 2\pi\ii\omega_\lambda\wedge)$ can recover 
the local system cohomology group $H^\bullet(\sfM(\scA), \scL)$.

Another interpretation of the Aomoto complex is 
related to the minimality of $\sfM(\scA)$ 
\cite{coh-int, dim-pap, ps-h, ran-mor}. 
In \cite{coh-int}, Cohen constructed a complex 
$(K^\bullet(\scA), \Delta^\bullet(\lambda))$ 
which computes $H^\bullet(\sfM(\scA), \scL)$ and 
the terms of this complex have Betti numbers 
as their dimensions, that is, satisfying the minimality: 
$\dim K^p=b_p(\sfM(\scA))$. 
When $\scL$ is trivial, the minimality implies 
all the coboundary maps vanish $\Delta^p(0)=0$. 
However the coboundary $\Delta(\lambda)$ of the 
minimal complex is difficult to compute for $\lambda\neq 0$. 
Cohen-Orlik \cite{coh-orl} determined the first order 
approximation of $\Delta(\lambda)$. 
They proved that the linearization 
of the minimal complex is chain equivalent to 
the Aomoto complex, 
$$
\left(K^\bullet(\scA), 
\left.
\frac{d}{dt}
\right|_{t=0}
\Delta^\bullet(t\lambda)
\right)
\cong 
(A^\bullet, 2\pi\ii\omega_\lambda\wedge). 
$$

The purpose of this paper is to study an ``integration'' 
of the Aomoto complex for obtaining 
the minimal complex for a real arrangement $\scA$. 
To do this, we introduce a basis of the 
Orlik-Solomon algebra $A^\bullet$, so called 
``chamber basis'', which is depending on a fixed 
generic flag with orientations. 
Then we have a matrix expression 
of the linear map 
$2\pi\ii\omega_\lambda\wedge:A^\bullet\rightarrow A^{\bullet +1}$. 
We will prove that the minimal complex can be recovered 
from these matrix entries. Roughly speaking, it is 
done just by replacing each matrix entry with 
its value of the hyperbolic sine function. 

The proof is based on the constructions in \cite{yos-min}. 
In the previous paper \cite{yos-min}, we 
explicitly constructed the attaching maps of 
cells arising from Lefschetz hyperplane section 
theorem for a real arrangement $\scA$. Moreover 
we also obtained a description of the 
minimal complex in terms of chambers. 
However the formula of the boundary map 
in \cite{yos-min} contains an 
integer $\deg(C', C)$ which is 
difficult to compute (see also Remark \ref{rem:deg}). 
The situation is changed 
in this paper. We will prove that the 
integer $\deg(C', C)$ appears as a ``structure constant'' 
of the Orlik-Solomon algebra with respect to the 
chamber basis. Moreover we also give an algorithm 
relating chamber basis and classical generator of 
the Orlik-Solomon algebra, which is more combinatorics friendly 
object than the original definition of $\deg(C',C)$.

The paper is organized as follows. In \S\ref{sec:pre}, 
we recall some basic facts on topology of $\sfM(\scA)$ and 
constructions from \cite{yos-min}. 
In particular, using a generic flag $\scF$, 
we divide the set of chambers into 
disjoint union $\ch(\scA)=\sqcup_{q=0}^\ell\ch_\scF^q(\scA)$ 
such that $|\ch_\scF^q(\scA)|=b_q(\sfM(\scA))$ 
and construct an isomorphism 
$\nu^q:\bZ[\ch_\scF^q(\scA)]\stackrel{\cong}{\longrightarrow}
H^q(\sfM(\scA), \bZ)$. 
This leads us to introduce the notion of chamber basis 
$\{\nu^q(C)| C\in\ch^q(\scA)\}$ of $A^q$. 
In \S\ref{sec:alg}, we will construct the inverse 
map 
$\xi^q=(\nu^q)^{-1}:H^q(\sfM,\bZ)
\rightarrow\bZ[\ch_\scF^q(\scA)]$ which enables us 
to express $\nu^q(C)\in H^q(\sfM(\scA))$ in terms of 
differential forms $\omega_H$. 
The wedge product 
$\omega_\lambda\wedge\nu^q(C)$ can be uniquely expressed as 
$\sum_{C'\in\ch^{q+1}}\Gamma_{C,C'}(\lambda)\nu^{q+1}(C')$ for 
some coefficients $\Gamma_{C,C'}(\lambda)\in\bC$. 
In \S\ref{sec:main}, 
two main results are stated and proved. 
First we assert that 
the coefficient $\Gamma_{C,C'}(\lambda)$ has a decomposition 
as a product of 
a linear form of weights $\lambda_H$ and an integer. 
Furthermore the linear factor of weights is explicitly 
described by using the notion of separating hyperplanes. 
The other integral factor is essentially equivalent to 
$\deg(C',C)$ mentioned above. 
The second result is recovering the minimal complex 
from these coefficients using the hyperbolic sine function. 
In the appendix, \S\ref{sec:app}, a generalized 
version of the linearization theorem for a minimal CW-complex 
is proved. 

\section{Preliminary}
\label{sec:pre}

\subsection{Basic constructions}
\label{sec:basic}

Let $V$ be an $\ell$-dimensional vector space. 
A finite set of affine hyperplanes 
$\scA=\{H_1, \ldots, H_n\}$ is called a {\it hyperplane arrangement}. 
For each hyperplane $H_i$ we fix a defining 
equation $\alpha_i$ such that $H_i=\alpha_i^{-1}(0)$. 
Let $L(\scA)$ be the set of nonempty intersections of 
elements of $\scA$. Define a partial order on $L(\scA)$ 
by $X\leq Y\Longleftrightarrow Y\subseteq X$ for $X, Y\in L(\scA)$. 
Note that this is reverse inclusion. 

Define a {\it rank function} on $L(\scA)$ by $r(X)=\codim X$. 
Write $L^p(\scA)=\{X\in L(\scA)|\ r(X)=p\}$. We call $\scA$ 
{\it essential} if $L^\ell(\scA)\neq\emptyset$. 

Let $\mu:L(\scA)\rightarrow \bZ$ be the {\it M\"obius function} 
of $L(\scA)$ defined by 
$$
\mu(X)=
\left\{
\begin{array}{ll}
1 &\mbox{ for }X=V\\
-\sum_{Y<X}\mu(Y), &\mbox{ for } X>V. 
\end{array}
\right.
$$
The {\it Poincar\'e polynomial} of $\scA$ is 
$\pi(\scA, t)=\sum\nolimits_{X\in L(\scA)}\mu(X)(-t)^{r(X)}$ and 
we also define numbers $b_i(\scA)$ by the formula 
$$
\pi(\scA, t)
=\sum\nolimits_{i=0}^\ell b_i(\scA)t^i. 
$$
We also define the {\it $\beta$-invariant} $\beta(\scA)$ by 
$$
\beta(\scA)=|\pi(\scA, -1)|, 
$$
if $\scA$ is an essential arrangement, the sign can be precisely 
enumerated as $\beta(\scA)=(-1)^\ell\pi(\scA, -1)$.

\subsection{Classical results}
\label{sec:classic}

Let $\scA$ be an arrangement in a real vector space 
$V_\bR$. Then the following relations between 
the set $\ch(\scA)$ of chambers and the complexified 
complement $\sfM(\scA)$ are known.

\begin{theorem}{\normalfont \cite{orl-sol, zas-face}}
\label{thm:classic}
\begin{itemize}
\item[{\normalfont (i) }] 
Let $\scA$ be an essential real $\ell$-arrangement. 
Then $|\ch(\scA)|=\pi(\scA, 1)$, and 
$|\bch(\scA)|=(-1)^\ell \pi(\scA, -1)=\beta(\scA)$, 
where $\bch(\scA)$ is the set of all bounded chambers.  
\item[{\normalfont (ii) }] 
Let $\scA$ be a complex arrangement. 
Then $b_i(\scA)$ is equal to the topological Betti number 
$b_i(\sfM(\scA))$, that is, 
$$
\Poin(\sfM(\scA), t)=\pi(\scA, t). 
$$
In particular, the absolute value of the 
topological Euler characteristic $|\chi(\sfM(\scA))|$ 
of the complement is equal to $\beta(\scA)$. 
\end{itemize}
\end{theorem}

\subsection{Generic flags and topology of $\sfM(\scA)$}
\label{sec:flag}

Let $\scA$ be an $\ell$-arrangement. 
A $q$-dimensional affine subspace $\scF^q\subset V$ is 
called {\it generic} or {\it transversal} to $\scA$ if 
$\dim \scF^q\cap X=q-r(X)$ for $X\in L(\scA)$. 
A {\it generic flag} 
$\scF$ is defined to be a complete flag (of affine subspaces) in $V$, 
$$
\scF:\ \emptyset=\scF^{-1}\subset\scF^0\subset\scF^1\subset\cdots 
\subset\scF^\ell=V,  
$$
where each $\scF^q$ is a generic $q$-dimensional affine subspace.

For a generic subspace $\scF^q$ we have an arrangement in $\scF^q$ 
$$
\scA\cap\scF^q:=\{H\cap\scF^q|\ H\in\scA\}. 
$$
The genericity provides an isomorphism of posets 
\begin{equation}
\label{eq:poset}
L(\scA\cap\scF^q)\cong L^{\leq q}(\scA):=\bigcup_{i\leq q}L^i(\scA). 
\end{equation}
In \cite{orl-sol} Orlik and Solomon gave a presentation 
of the cohomology ring $H^*(\sfM(\scA), \bZ)$ in terms 
of the poset $L(\scA)$ for a complex arrangement $\scA$. 
The next proposition follows from (\ref{eq:poset}). 
\begin{proposition}
\label{prop:trunc}
Let $\scA$ be a complex arrangement and $\scF^q$ a 
$q$-dimensional generic subspace. Then the natural inclusion 
$i:\sfM(\scA)\cap\scF^q\hookrightarrow \sfM(\scA)$ 
induces isomorphisms 
$$
i_k: H_k(\sfM(\scA)\cap\scF^q, \bZ)\stackrel{\cong}{\longrightarrow}
H_k(\sfM(\scA), \bZ), 
$$
for $k=0, 1, \ldots, q$. 
\end{proposition}
In particular, the Poincar\'e polynomial of $\scA\cap\scF^q$ 
is given by 
\begin{equation}
\label{eq:trun}
\pi(\scA\cap\scF^q, t)=\pi(\scA, t)^{\leq q}, 
\end{equation}
where $(\sum_{i\geq 0}a_i t^i)^{\leq q}=\sum_{i=0}^q a_i t^i$ is the 
truncated polynomial. From these formulas and 
Theorem \ref{thm:classic}, we have the following 
proposition. (For the proof see \cite[Prop. 2.3.2]{yos-min} 
for example.) 

\begin{proposition}
\label{prop:flag}
Let $\scA$ be a real $\ell$-arrangement and $\scF$ a generic flag. 
Define 
$$
\sfch_\scF^q(\scA)=\{C\in\sfch(\scA)|\ C\cap\scF^q\neq\emptyset 
\mbox{ {\normalfont and} }C\cap\scF^{q-1}=\emptyset \}, 
$$
for each $q=0, 1, \ldots, \ell$. 
Then $|\ch_\scF^q(\scA)|=b_q(\sfM(\scA)) $. 
\end{proposition}
In particular, the number of chambers which does not 
intersect with a generic hyperplane $\scF^{\ell-1}$ 
satisfies 
$$
|\ch_\scF^\ell(\scA)|=b_\ell(\sfM(\scA)).
$$
This formula has a topological meaning. 
Let us recall the construction in \cite{yos-min} 
briefly. 
First to fix an orientation, we fix a basis $(v_1, \ldots, v_\ell)$ 
of $V$ such that 
$$
\scF^q=\scF^0+\sum_{i=1}^q\bR v_i. 
$$
The orientation of $\scF^q$ is determined by 
the ordered basis $(v_1, \ldots ,v_q)$. And also 
define positive and negative half spaces, 
$\scF_+^q$ and $\scF_-^q$, by 
\begin{eqnarray*}
\scF_+^q&=&\scF^{q-1}+\bR_{>0}v_q\\
\scF_-^q&=&\scF^{q-1}+\bR_{<0}v_q, 
\end{eqnarray*}
respectively. 

\begin{definition}
The map $\sign: \ch_\scF^q(\scA)\rightarrow \{\pm 1\}$ 
is defined by 
$$
\sign(C)=
\left\{
\begin{array}{cl}
1& \mbox{ if } \scF^q\cap C\subset \scF_+^q\\
-1& \mbox{ if } \scF^q\cap C\subset \scF_-^q. 
\end{array}
\right.
$$
\end{definition}

Let $\scF^q\subset\bR^\ell$. 
Denote by $\scF_\bC^q=\scF^q\otimes\bC$ 
the complexification of $\scF^q$ and put 
$\sfM^q :=\scF_\bC^q\cap\sfM(\scA)$. Note that
$\sfM^\ell=\sfM(\scA)$. 
We fix the orientation of $\scF_\bC^q$ by the 
ordered basis 
$$
(v_1, \ldots, v_q, \sqrt{-1}v_1, \ldots,\sqrt{-1}v_q).
$$ 
Note that this orientation is different by 
$(-1)^{\frac{q(q-1)}{2}}$ from the canonical orientation 
of a complex vector space. 

For each $C\in\ch_\scF^\ell(\scA)$, we can explicitly 
construct a continuous map 
$$
\sigma_C:
(D^\ell, \partial D^\ell)
\longrightarrow 
(\sfM^\ell, \sfM^{\ell-1}), 
$$
such that (\cite[\S 5.2]{yos-min})
\begin{equation}
\label{eq:cell}
\begin{array}{l}
\mbox{(Transversality) }\sigma_C(0)\in C, 
\sigma_C(D^\ell)\pitchfork C=\{\sigma_C(0)\} \mbox{ and }\\
\\
\mbox{(Non-intersecting) }\sigma_C(D^\ell)\cap C'=\emptyset \mbox{ for }
C'\in\ch_\scF^\ell(\scA)\backslash\{C\}. 
\end{array}
\end{equation}
These properties guarantee the following homotopy 
equivalence (\cite[4.3.1]{yos-min}): 
\begin{equation}
\label{eq:attach}
\sfM^\ell\simeq\sfM^{\ell-1}\cup_{(\partial\sigma_C)}
\left(
\bigsqcup_{C\in\ch_\scF^\ell(\scA)}D^\ell
\right), 
\end{equation}
where the right hand side is obtained by 
attaching $\ell$-dimensional disks to $\sfM^{\ell-1}$ 
along 
$\partial\sigma_C:\partial D^\ell\rightarrow \sfM^{\ell-1}$ 
for ${C\in\ch_\scF^\ell(\scA)}$. 
Since the natural inclusion $\sfM^{\ell-1}\hookrightarrow\sfM^\ell$ 
induces 
$H_{\ell-1}(\sfM^{\ell-1}, \bZ)\cong H_{\ell-1}(\sfM^{\ell},\bZ)$ 
and 
$H_\ell(\sfM^\ell, \bZ)\cong H_\ell(\sfM^{\ell},\sfM^{\ell-1}, \bZ)$, 
$\sigma_C(D^\ell)$ 
can be considered as an element of 
$H_\ell(\sfM^{\ell}, \bZ)$. 
Furthermore, 
$\{\sigma_C(D^\ell)\}_{C\in\ch_\scF^\ell}$ form a basis of 
$H_\ell(\sfM^\ell, \bZ)$. 
We choose an orientation of $\sigma_C$ so 
that the intersection number satisfies 
$$
[C]\cdot [\sigma_C]=1. 
$$
Then chambers 
$\{[C]\}_{C\in\ch_\scF^\ell}$ form the dual basis of 
locally finite homology group 
$H_\ell^{lf}(\sfM^\ell, \bZ)$, which is isomorphic to 
$H^\ell(\sfM^\ell, \bZ)$. Thus we have the 
isomorphism 
$$
\begin{array}{ccccccc}
\bZ[\ch_\scF^q(\scA)]&
\stackrel{\cong}{\longrightarrow}&
H_q^{lf}(\sfM^q, \bZ)&
\stackrel{\cong}{\longrightarrow}&
H^q(\sfM^q, \bZ)&
\stackrel{\cong}{\longrightarrow}&
H^q(\sfM^\ell, \bZ)\\
&&&&&&\\
C&&&&&&\nu^q(C). 
\end{array}
$$
Denote by $\nu^q$ the composite map 
\begin{equation}
\label{eq:isom}
\nu^q: \bZ[\ch_\scF^q(\scA)]
\stackrel{\cong}{\longrightarrow}
H^q(\sfM^\ell, \bZ). 
\end{equation}

\begin{definition}
The set $\{\nu^q(C)|C\in\ch_\scF^q(\scA)\}$ is 
called the chamber basis of $H^q(\sfM, \bZ)$ with 
respect to a flag $\scF$. 
\end{definition}

\begin{remark}
In \cite{var-gel}, 
Varchenko and Gel'fand constructed a filtration 
$0\subset P^0\subset\cdots\subset P^\ell=\bZ[\ch(\scA)]$ 
and an isomorphism 
$P^q/P^{q-1}\cong H_{2\ell-q}^{lf}(\sfM(\scA), \bZ)$. 
Our subspace $\bZ[\ch_\scF^\ell(\scA)]$ gives a section 
of the quotient map $P^\ell\rightarrow P^\ell/P^{\ell-1}$.  
Using generic section by $\scF^q$, we can 
construct an isomorphism 
$\bZ[\ch_\scF^q(\scA)]\cong P^q/P^{q-1}$. 
The map $\nu^q$ is equivalent to 
Varchenko-Gel'fand's isomorphism under the identification 
$H_{2\ell-q}^{lf}(\sfM(\scA), \bZ)\cong H^q(\sfM(\scA), \bZ)$ 
up to sign. 
\end{remark}



\section{An algorithm relating chamber and 
differential forms}
\label{sec:alg}

From the result of Brieskorn, Orlik and Solomon, 
the cohomology ring $A^*=H^*(\sfM(\scA),\bZ)$ is 
generated by 
$\omega_i=\frac{1}{2\pi\ii}\frac{d\alpha_i}{\alpha_i}$ 
($i=1, \ldots, n$). 
In this section, we express $\nu^q(C)$ in terms of generators 
$\omega_H$. 

Let $I=\{i_1, \ldots, i_\ell\}\subset\{1, \ldots, n\}$ be 
an ordered subset of $\ell$ indices, 
$\scA(I):=\{H_{i_1}, \ldots, H_{i_\ell}\}$ be 
a subarrangement consists of $\ell$ hyperplanes. 
Suppose $H_{i_1}, \ldots, H_{i_\ell}$ are independent, that is, 
$d\alpha_{i_1}\wedge\cdots\wedge d\alpha_{i_\ell}\neq 0$. 
Obviously $\ch(\scA(I))$ consists of 
$2^\ell$ chambers and there exists a unique 
chamber $C_0(I)
\in\ch(\scA(I))$ with 
$C_0(I)\cap \scF^{\ell-1}=\emptyset$. 
Choose a normal vector $w_{i_k}\perp H_{i_k}$ for 
each $H_{i_k}$ such that $C_0(I)$ is contained in 
the half space $H_{i_k}+\bR_{>0}\cdot w_{i_k}$. 

\begin{definition}
For an ordered $\ell$-tuple 
$I=(i_1, \ldots, i_\ell)\subset\{1, \ldots, n\}$, define 
$\varepsilon(I)$ by 
$$
\varepsilon(I)=
\left\{
\begin{array}{ll}
0& \mbox{if $H_{i_1}, \ldots, H_{i_\ell}$ are dependent,} \\
1& \mbox{if $(w_{i_1}, \ldots, w_{i_\ell})$ is a positive basis,} \\
-1&\mbox{if $(w_{i_1}, \ldots, w_{i_\ell})$ is a negative basis.}  
\end{array}
\right.
$$
\end{definition}
For an ordered $\ell$-tuple 
$I=(i_1, \ldots, i_\ell)$, denote 
$\omega_I=\omega_{i_1}\wedge\cdots\wedge\omega_{i_\ell}$. 
Let us define as follows. 
$$
\xi^\ell(\omega_I)=\varepsilon(I)[C_0(I)]\in
\bZ[\ch_\scF^\ell(\scA)]. 
$$ 

\begin{theorem}
\label{thm:xi}
The map $\xi^\ell$ induces the isomorphism 
$$
\xi^\ell:H^\ell(\sfM(\scA), \bZ)
\stackrel{\cong}{\longrightarrow}
\bZ[\ch_\scF^\ell(\scA)], 
$$
and $\xi^\ell=(\nu^\ell)^{-1}$. 
\end{theorem}

\noindent
{\bf Proof.} We prove $\nu^\ell(\xi^\ell(\omega_I))=\omega_I$. 
It is enough to show 
$$
\int_{\sigma_C}\omega_I=
\varepsilon(I)[C_0(I)]\cdot\sigma_C, 
$$
for $C\in\ch_\scF^\ell(\scA)$. 
Note that the inclusion 
$\iota:\sfM(\scA)\hookrightarrow\sfM(\scA(I))$ 
induces a surjective map 
$\iota_*:H_\ell(\sfM(\scA),\bZ)\rightarrow 
H_\ell(\sfM(\scA(I)),\bZ)$. 
Since $\sfM(\scA(I))$ is homeomorphic to 
$(\bC^*)^\ell$, the top homology 
$H_\ell(\sfM(\scA(I)), \bZ)$ is generated by 
the cycle 
$$
T:=\{(\alpha_{i_1}, \ldots, \alpha_{i_\ell}); 
|\alpha_{i_1}|=|\alpha_{i_2}|=\cdots =|\alpha_{i_\ell}|=1\}. 
$$
Fix the orientation of $T$ 
so that $[C_0(I)]\cdot T=1$. 
Then since $[C_0]$ is the dual basis to $T$, 
from (\ref{eq:cell}), we have, 
$$
\iota_*(\sigma_C)=
\left\{
\begin{array}{ll}
T& \mbox{if } C\subseteq C_0(I),\\
0& \mbox{if } C\cap C_0(I)=\emptyset, 
\end{array}
\right.
$$
for $C\in\ch_\scF^\ell(\scA)$. 
Now we have $\int_{\sigma_C}\omega_I=\int_{\iota_*(\sigma_C)}\omega_I$ 
and 
$\int_T\omega_I=\varepsilon(I)$ completes the proof. 
\owari

Similarly, we can define 
$\xi^q:H^q(\sfM(\scA), \bZ)
\stackrel{\cong}{\longrightarrow}
\bZ[\ch_\scF^q(\scA)]$ by using 
$H^q(\sfM(\scA), \bZ)\cong H^q(\sfM(\scA)\cap\scF^q, \bZ)$ 
for $0\leq q\leq \ell-1$. 
Theorem \ref{thm:xi} enables us to express 
chamber basis $\nu^q(C)$ in terms of generators 
$\omega_H$. 

\begin{example}
Let $\scA=\{H_1, \ldots, H_4\}$ be the 
arrangement of $4$-lines as in Figure \ref{fig:ex} with 
flag $\scF^\bullet$ defined by $v_1, v_2$. Then 
$\ch_\scF^0(\scA)=\{A\}$, 
$\ch_\scF^1(\scA)=\{B_1, B_2, B_3, B_4\}$, and 
$\ch_\scF^2(\scA)=\{C_1, C_2, C_3, C_4, C_5\}$. 
From Theorem \ref{thm:xi}, $\xi(1)=[A], \nu(A)=1$ and 
$$
\begin{array}{rclcrcl}
\xi(\omega_1)&=&-[B_1]&&
\nu(B_1)&=&-\omega_1\\
\xi(\omega_2)&=&[B_2]+[B_3]+[B_4]&&
\nu(B_2)&=&\omega_2-\omega_3\\
\xi(\omega_3)&=&[B_3]+[B_4]&&
\nu(B_3)&=&\omega_3-\omega_4\\
\xi(\omega_4)&=&[B_4]&&
\nu(B_4)&=&\omega_4
\end{array}
$$

$$
\begin{array}{rclcrcl}
\xi(\omega_{12})&=&-[C_1]-[C_3]&&
\nu(C_1)&=&-\omega_{12}+\omega_{14}-\omega_{24}\\
\xi(\omega_{13})&=&-[C_1]-[C_2]-[C_3]-[C_4]&&
\nu(C_2)&=&\omega_{12}-\omega_{13}+\omega_{24}-\omega_{34}\\
\xi(\omega_{14})&=&-[C_3]-[C_4]-[C_5]&&
\nu(C_3)&=&-\omega_{14}+\omega_{24}\\
\xi(\omega_{24})&=&-[C_4]-[C_5]&&
\nu(C_4)&=&-\omega_{24}+\omega_{34}\\
\xi(\omega_{34})&=&-[C_5]&&
\nu(C_5)&=&-\omega_{34}. 
\end{array}
$$

$$
\omega_\lambda\wedge
\nu(A)=
-\lambda_1\nu(B_1)+\lambda_2\nu(B_2)
+\lambda_{23}\nu(B_3)+\lambda_{234}\nu(B_4). 
$$

$$
\begin{array}{crrrrr}
\omega_\lambda\wedge\nu(B_1)=&
-\lambda_{23}\nu(C_1)&
-\lambda_{3}\nu(C_2)&
-\lambda_{234}\nu(C_3)&
-\lambda_{34}\nu(C_4)&
-\lambda_{4}\nu(C_5)\\
\omega_\lambda\wedge\nu(B_2)=&
&
+\lambda_{123}\nu(C_2)&
&
+\lambda_{1234}\nu(C_4)&
\\
\omega_\lambda\wedge\nu(B_3)=&
-\lambda_{1}\nu(C_1)&
-\lambda_{12}\nu(C_2)&
&
&
+\lambda_{1234}\nu(C_5)\\
\omega_\lambda\wedge\nu(B_4)=&
&
&
-\lambda_{1}\nu(C_3)&
-\lambda_{12}\nu(C_4)&
-\lambda_{123}\nu(C_5)
\end{array}
$$

\begin{figure}[htbp]
\begin{picture}(100,150)(20,0)
\thicklines

\put(200,0){\line(0,1){150}}
\put(80,0){\line(2,1){250}}
\put(248,0){\line(-4,5){120}}
\put(70,150){\line(5,-2){350}}

\multiput(50,20)(5,0){73}{\circle*{2}}
\put(170,20){\circle*{5}}
\put(170,20){\vector(1,0){20}}
\put(170,20){\vector(0,1){20}}
\put(158,6){$\scF^0$}
\put(33,16){$\scF^1$}
\put(185,10){$v_1$}
\put(172,35){$v_2$}

\put(80,-12){$H_1$}
\put(195,-12){$H_2$}
\put(250,-12){$H_3$}
\put(410,0){$H_4$}

\put(130,0){$A$}
\put(70,50){$B_1$}
\put(215,0){$B_2$}
\put(280,40){$B_3$}
\put(370,60){$B_4$}

\put(210,75){$C_1$}
\put(185,85){$C_2$}
\put(240,110){$C_3$}
\put(175,120){$C_4$}
\put(110,140){$C_5$}

\end{picture}
     \caption{Example}\label{fig:ex}
\end{figure}
\end{example}


\section{Aomoto complex via chamber-basis}
\label{sec:main}

\subsection{Main result}

Let $(\lambda_1, \ldots, \lambda_n)\in\bC^n$ and 
put $\omega_\lambda=\sum_{i=1}^n\lambda_i\omega_i$. 
Since $\omega_\lambda\wedge\omega_\lambda=0$, we have 
a cochain complex $(A^\bullet, 2\pi\ii\omega_\lambda\wedge)$, 
which is called the Aomoto complex. We shall 
study this complex using the chamber basis 
$\{\nu^q(C)\}_{C\in\ch_\scF^q(\scA)}$ of $A^*$. 
For a chamber $C\in\ch_\scF^q(\scA)$, 
$\omega_\lambda\wedge\nu^q(C)$ 
is uniquely expressed as 
$$
\omega_\lambda\wedge\nu^q(C)=
\sum_{C'\in\ch_\scF^{q+1}(\scA)}\Gamma_{C,C'}(\lambda)\cdot
\nu^{q+1}(C'), 
$$
for some complex numbers 
$\Gamma_{C,C'}(\lambda)\in\bC$. We may consider the 
coefficients $\{\Gamma_{C,C'}(\lambda)\}$ as structure constants 
of the cohomology ring with respect to the chamber 
basis. 

Let $\scL_\lambda$ be a rank one local system on $\sfM(\scA)$ 
determined by 
monodromies $q_i=e^{2\pi\ii \lambda_i}\in\bC^*$ 
around the hyperplane $H_i$. 
For given two chambers $C, C'\in\ch(\scA)$, denote 
$S(C, C')$ the set 
$$
S(C, C')=\{H\in\scA\ |\ \mbox{$H$ separates $C$ and $C'$}\}, 
$$
of hyperplanes separating $C$ and $C'$, and 
$\lambda_{S(C, C')}=\sum_{H\in S(C, C')}\lambda_H$. 

The main result is the following. 

\begin{theorem}
(a) The coefficient $\Gamma_{C,C'}(\lambda)$ has the following 
decomposition. 
$$
\Gamma_{C,C'}(\lambda)=N_{C,C'}\cdot\lambda_{S(C,C')}, 
$$
where $N_{C,C'}\in\bZ$. \\
(b) 
Let us define a linear map 
$\wt{\nabla}_\lambda:\bC[\ch_\scF^q(\scA)]\rightarrow
\bC[\ch_\scF^{q+1}(\scA)]$ by 
$$
\wt{\nabla}_\lambda([C])=
-\sum_{C'\in\ch_\scF^{q+1}(\scA)}2\cdot N_{C,C'}\cdot
\sinh\left(
\pi\ii{\lambda_{S(C,C')}}\right)
[C']. 
$$
Then $(\ch_\scF^\bullet(\scA), \wt{\nabla}_\lambda)$ is a cochain 
complex and 
$$
H^p(\ch_\scF^\bullet(\scA), \wt{\nabla}_\lambda)\cong 
H^p(\sfM(\scA), \scL_\lambda). 
$$
\end{theorem}
Thus $(\ch_\scF^\bullet(\scA), \wt{\nabla}_\lambda)$ is 
the minimal complex which is obtained as a 
integration of the Aomoto complex.

\subsection{Proof}

First we recall some more notation 
from \cite{yos-min}. 
We defined the degree map (\cite[\S 6.3]{yos-min}) 
$$
\deg:\ch_\scF^{p+1}(\scA)\times\ch_\scF^{p}(\scA)
\longrightarrow \bZ. 
$$
Furthermore the notion of $\deg(C',C)$ enables us 
to express twisted cellular coboundary map 
$\nabla_{\scL_\lambda}:\bC[\ch_\scF^p(\scA)]
\rightarrow \bC[\ch_\scF^{p+1}(\scA)]$ as follows 
(\cite[6.4.1]{yos-min}): 
$$
\nabla_{\scL_\lambda}(C)=
-\sum_{C'\in\ch_\scF^{p+1}(\scA)}\sign(C')\deg(C',C)\cdot 
2\sinh\left(\pi\ii\lambda_{S(C',C)}\right)[C']. 
$$
In particular, 
$(\bC[\ch_\scF^\bullet], \nabla_{\scL_\lambda})$ is a cochain complex, 
and the cohomology group is isomorphic to cohomology 
with coefficients in ${\scL_\lambda}$: 
$H^p(\bC[\ch_\scF^\bullet], \nabla_{\scL_\lambda})\cong 
H^p(\sfM(\scA), \scL_\lambda)$. 

We now apply the linearization theorem by 
Cohen-Orlik \cite{coh-orl}. 
Consider the local system $\scL_{t\lambda}$ with $t\in\bC$. 
Since 
\begin{equation}
\label{eq:sinh}
\left.
\frac{d}{dt}
\right|_{t=0}
2\sinh\left(\pi\ii t\lambda_{S(C',C)}\right)
=
2\pi\ii\lambda_{S(C,C')}, 
\end{equation}
we have 
$$
\left.
\frac{d}{dt}
\right|_{t=0}
\nabla_{\scL_{t\lambda}}(C)=
-2\pi\ii\sum_{C'\in\ch_\scF^{p+1}}\sign(C')\deg(C',C)
\lambda_{S(C,C')}[C']. 
$$
From the construction in \cite[6.4]{yos-min}, 
$[C]$ can be identified with $\nu(C)$ here. 
Thus we have 
$\Gamma_{CC'}(\lambda)=-\sign(C')\deg(C',C)\lambda_{S(C,C')}$. 
The map $\wt{\nabla}_\lambda$ of (b) is clearly equivalent to 
${\nabla}_{\scL_\lambda}$. 
\owari

\begin{remark}
\label{rem:deg}
In \cite{yos-min}, a minimal CW-decomposition 
such that each $k$-cell is labeled by a chamber 
$C\in\ch_\scF^k(\scA)$ is constructed. From the 
minimality, the incidence numbers vanish $[C':C]=0$. 
On the other hand, minimal CW-decomposition of $\sfM(\scA)$ 
induces a $\bZ^{b_1}$-equivariant CW-decomposition of 
the homology covering $\wt{\sfM}$. 
The degree above can be considered as the incidence 
number at the level of homology covering. 
\end{remark}


\section{Appendix: Linearization theorem 
for minimal CW-complex}

\label{sec:app}

In this section, we give a proof of 
the linearization theorem by Cohen-Orlik \cite{coh-orl} 
in a generalized setting. 

Let $X$ be a connected minimal CW-complex, that is, 
a finite CW-complex with exactly as many $k$-cells 
as the $k$-th Betti number, for all $k$. 
Denote by $\scS_k$ the set of $k$-cells and 
by $X_k$ the $k$-skeleton of $X$. 
Suppose $|\scS_1|=n$, then 
$H_1(X, \bZ)\cong\bZ^n$ and $H^1(X, \bZ)=\bZ^n$. 
An example of a minimal CW-complex is 
the $n$-torus $T^n=(S^1)^{\times n}$. 
The $n$-torus $T^n$ admits the canonical 
minimal cell decomposition as follows. 
Let $e_1, \ldots, e_n$ be the standard basis 
of $\bR^n$. The torus $T^n$ can be identified with 
the quotient space 
$\left(\bigoplus_{i=1}^n\bR e_i\right)/
\left(\bigoplus_{i=1}^n\bZ e_i\right)$. 
For any subset 
$\Phi=\{p_1, \ldots, p_k\}\subset[1,n]:=\{1, \ldots, n\}$, 
denote by $K_\Phi$ the $k$-cube 
$$
K_\Phi=
\{t_1e_{p_1}+\cdots + t_ke_{p_k}| 0\leq t_i\leq 1, i=1, \ldots, k\}
\subset\bR^n 
$$
and by $e_\Phi=p(K_\Phi)$ the image of the cube 
by the quotient map 
$p:\bR^n\rightarrow \bR^n/\bZ^n$. This gives a 
cell decomposition of $T^n=\bigcup_{\Phi}e_\Phi$. 
Note that $\dim e_\Phi=|\Phi|$. 
Obviously the quotient map $p$ gives the universal covering 
of $T^n$. Let us denote by $x_i$ 
the multiplicative generator of the deck transformation group 
corresponding to $e_i$. The deck transformation group 
is identified with the multiplicative group of Laurent monomials 
$\{x^\alpha=x_1^{\alpha_1}x_2^{\alpha_2}\cdots x_n^{\alpha_n}| 
\alpha=(\alpha_1, \ldots, \alpha_n)\in\bZ^n\}$. 
The covering space $\wt{T}^n\cong\bR^n$ has the following 
cell decomposition, 
$$
\wt{T}^n=\bigcup_{\Phi\subset [1,n]}\bigcup_{\alpha\in H_1}x^\alpha\cdot K_\Phi. 
$$

Since $T^n$ is the $K(\bZ^n, 1)$-space, 
the abelianization map 
$\pi_1(X)\rightarrow H_1(X, \bZ)\cong \bZ^n$ 
determines, uniquely up to homotopy, 
a continuous map $f:X\rightarrow T^n$ 
such that $f_*:H_1(X, \bZ)\stackrel{\cong}{\longrightarrow}
H_1(T^n, \bZ)$. By cellular approximation theorem, we may 
assume that $f$ is cellular, i.e., 
preserving skeletons $f(X_k)\subset (T^n)_k$. 
A $k$-cell $\sigma\in\scS_k$ is expressed as 
a characteristic map 
$$
\sigma:(D^{k}, \partial D^k)\longrightarrow (X_k, X_{k-1}), 
$$
from the $k$-disk to the $k$-skeleton. 
For simplicity, we assume the base point 
$p_\sigma\in D^k$ 
is mapped to $X_0$ by $\sigma$. 

As is the case of $T^n$, 
the $H_1(X, \bZ)$-covering $\widetilde{X}$ of $X$ 
has the structure of $\bZ^n$-equivariant CW-complex. 
Indeed from the minimality, 
fixing a base point $p_0\in\widetilde{X}$ over 
$X_0$, each cell $\sigma:D^k\rightarrow X$ can be lifted uniquely 
$\widetilde{\sigma}:D^k\rightarrow\widetilde{X}$ such that 
$\widetilde{\sigma}(p_\sigma)=p_0$. Then $\widetilde{X}$ 
is decomposed as 
$$
\widetilde{X}=
\bigcup_{\sigma\in\scS}
\bigcup_{\alpha\in H_1}x^\alpha\cdot
\widetilde{\sigma},  
$$
where $\scS=\bigcup_i\scS_i$ is the set of all cells. 
The boundary of the cellular chain complex can be expressed 
as 
\begin{equation}
\label{eq:cover}
\partial(x^\alpha\cdot\widetilde{\sigma})=
\sum_{\tau\in\scS_{k-1}}\sum_{\beta\in H_1}
[\widetilde{\sigma} :  (x^\beta\cdot\widetilde{\tau})]
x^{\alpha+\beta}\cdot\widetilde{\tau}
\in H_k(\wt{X}_k, \wt{X}_{k-1}; \bZ), 
\end{equation}
where $[\widetilde{\sigma} :  (x^\beta\cdot\widetilde{\tau})]\in\bZ$ 
is the incidence number.

Since $f$ induces the isomorphism 
$H_1(X, \bZ)\cong H_1(T^n, \bZ)$, any complex rank one 
local system on $X$ can be obtained as a pull back 
by $f$. 
Recall that a local system on $T^n$ is determined by 
a homomorphism 
$\rho:H_1(T^n, \bZ)\rightarrow\bC^*$. 
Let $\lambda=(\lambda_1, \ldots, \lambda_n)\in\bC^n$, 
and define a local system $\scL_\lambda$ on $T^n$ by 
$\rho(e_{\{i\}})=q_i=e^{2\pi\ii\lambda_i}\in\bC^*$. 
We also denote $\rho(x^\alpha)=q^\alpha=\prod_{i=1}^nq_i^{\alpha_i}$. 
Now we describe the boundary map of a $\scL_\lambda$-coefficients 
cellular complex. Recall that (\cite[VI.2]{whi-ele}) 
a chain with coefficients in a local system 
$\scL_{\lambda}$ is 
a pair of a cell $\sigma$ and a section 
$c(\sigma)\in(\scL_{\lambda})_{\sigma(p_\sigma)}$ at the base point. 
Let us fix a trivialization 
$(\scL_{\lambda})_{X_0}\cong\bC$ at the base point. We identify the 
cell $\sigma$ with a chain 
$[\sigma]$ with $\scL_{\lambda}$-coefficient which takes 
the value $c=1\in\bC\cong(\scL_{\lambda})_{X_0}$. 
Then, for a $k$-cell $\sigma\in\scS_k$, the boundary 
with $\scL_{\lambda}$-coefficients is expressed as 
\begin{equation}
\label{eq:bdry}
\partial_{\scL_{\lambda}}[\sigma]=
\sum_{\tau\in\scS_{k-1}}\sum_{\alpha\in H_1}
[\widetilde{\sigma}:(x^\alpha\cdot{\tau})]
q^\alpha[\tau]
\in H_k({X}_k, {X}_{k-1}; \bC). 
\end{equation}

Let us denote the dual basis of $e_{\{i\}}$ by 
$dt_i\in H^1(T^n,\bZ)$ and define 
$$
\omega_\lambda=2\pi\ii\sum_{i=1}^n\lambda_idt_i\in H^1(T^n, \bC).
$$ 
Clearly the minimality implies 
$\lim_{t\rightarrow 0}\partial_{\scL_{t\lambda}}[\sigma]=0$. 
By differentiating (\ref{eq:bdry}), we have the following 
lemma. 

\begin{lemma}
\label{lem:diff}
Let $X$ be a minimal CW-complex. 
Let $t\in\bC$ and 
$\sigma\in\scS_k$ be a $k$-cell of $X$. 
\begin{equation}
\label{eq:diff}
\left.
\frac{d}{dt}
\right|_{t=0}
\partial_{\scL_{t\lambda}}[\sigma]=
\sum_{\tau\in\scS_{k-1}}\sum_{\alpha\in H_1}
[\widetilde{\sigma}:(x^\alpha\cdot{\tau})]
\cdot 
\langle\omega_{\lambda}, \alpha\rangle
\cdot[\tau]. 
\end{equation}
\end{lemma}
The next result asserts that, after 
push forward by $f_*$,  it is obtained 
by the cap product on the torus. 
\begin{theorem}
\begin{equation}
\label{eq:bibun}
f_*
\left(
\left.
\frac{d}{dt}
\right|_{t=0}
\partial_{\scL_{t\lambda}}[\sigma]
\right)
=
\omega_\lambda\cap f_*[\sigma]. 
\end{equation}
\end{theorem}

\noindent
{\bf Proof.}  
We prove this in two steps. The first step is 
commuting the operators $f_*$ and $\ddt\partial_{\scL_{t\lambda}}$. 
This is essentially done by the fact that a cellular 
map induces a chain map between cellular chain complexes. 
Therefore it is enough to show that 
\begin{equation}
\label{eq:goal}
\ddt\partial_{\scL_{t\lambda}}[f\circ\sigma]=
\omega_\lambda\cap [f\circ\sigma]. 
\end{equation}
Let us denote by $\eta$ the left hand side of 
this formula. 
Note that the continuous map $f:X\rightarrow T^n$ can also be 
lifted to $\wt{f}:\wt{X}\rightarrow \wt{T}^n$. 
From Lemma \ref{lem:diff}, we have 
$$
\eta =
\sum_{|\Phi|=k-1}
\sum_{\alpha\in H_1}[\widetilde{f}\circ\widetilde{\sigma}: 
(x^\alpha\cdot e_\Phi)]\langle\omega_\lambda, \alpha\rangle [e_\Phi]. 
$$
Let $\Phi_0\subset[1,n]$ with $|\Phi_0|=k-1$. 
Then (\ref{eq:goal}) is equivalent to 
\begin{equation}
\label{eq:int}
\int_\eta dt_{\Phi_0}=
\int_{[\wt{f}\circ\wt{\sigma}]}\omega_\lambda\wedge dt_{\Phi_0}, 
\end{equation}
for all $\Phi_0$. 
The left hand side of (\ref{eq:int}) is equal to 
$$
\sum_\alpha[\wt{f}\circ\wt{\sigma}: 
(x^\alpha\cdot e_{\Phi_0})]\langle\omega_\lambda, \alpha\rangle.  
$$
We compute the right hand side of (\ref{eq:int}) using 
$\omega_\lambda=d(2\pi\ii\sum\lambda_it_i)$, 
by Stokes theorem, 
\begin{eqnarray*}
\int_{[\wt{f}\circ\wt{\sigma}]}\omega_\lambda\wedge dt_{\Phi_0}&=&
\sum_{i=1}^n\int_{\partial[\wt{f}\circ\wt{\sigma}]}2\pi\ii\lambda_it_i dt_{\Phi_0}\\
&=&
2\pi\ii\sum_{i=1}^n\sum_{|\Phi|=k-1}\sum_{\alpha\in H_1}
\int_{x^\alpha e_{\Phi}}[(\wt{f}\circ\wt{\sigma}): (x^\alpha e_{\Phi})]
\lambda_it_i dt_{\Phi_0}\\
&=&
2\pi\ii\sum_{i=1}^n\sum_{\alpha\in H_1}
\int_{e_{\Phi_0}}[(\wt{f}\circ\wt{\sigma}): (x^\alpha e_{\Phi_0})]
\lambda_i\alpha_i dt_{\Phi_0}\\
&=&
\sum_{\alpha\in H_1}
[(\wt{f}\circ\wt{\sigma}): (x^\alpha e_{\Phi_0})]
\langle\omega_\lambda, \alpha\rangle. 
\end{eqnarray*}
This completes the proof of (\ref{eq:int}). \owari

\begin{corollary}
If $X$ satisfies that 
$f_*:H_k(X, \bZ)\rightarrow H_k(T^n, \bZ)$ injective, 
then 
\begin{equation}
\label{eq:lin}
\left.
\frac{d}{dt}
\right|_{t=0}
\partial_{\scL_{t\lambda}}[\sigma]
=
f^*(\omega_\lambda)\cap \sigma, 
\end{equation}
for every $k$-cell $\sigma\in\scS_k$. 
\end{corollary}

\begin{remark}
If $X$ is the complement to a hyperplane arrangement, 
then $f_*$ is injective for every $k$. Thus we have 
the linearization theorem by \cite{coh-orl}. 
The author does not know whether the formula (\ref{eq:lin}) 
holds for any minimal CW-complex $X$. 
\end{remark}

\medskip

\noindent
{\bf Acknowledgement.} 
The idea considering the chamber basis of 
the Orlik-Solomon algebra 
was inspired from correspondences with 
Professor H. Terao and Dr. T. Abe. 
The author thanks to them. 
The author is supported by JSPS Postdoctoral Fellowship 
Research Abroad. He wishes to express his gratitude 
to JSPS for the support and 
to The Abdus Salam ICTP for wonderful working enviroment 
and hospitality.

\noindent
Masahiko Yoshinaga

The Abdus Salam International Centre for Theoretical Physics,
Strada Costiera 11,
Trieste 34014,
Italy. email: myoshina@ictp.it


\begin{thebibliography}{99}

\bibitem[C]{coh-int}
D. Cohen, 
Cohomology and intersection cohomology of complex hyperplane arrangements. 
Adv. Math. {\bf 97} (1993), no. 2, 231--266.

\bibitem[CO]{coh-orl}
D. Cohen, P. Orlik, 
Arrangements and local systems. 
Math. Res. Lett. {\bf 7} (2000), no. 2-3, 299--316.

\bibitem[CS]{coh-suc}
D. Cohen, A. Suciu, 
Characteristic varieties of arrangements.  
Math. Proc. Cambridge Philos. Soc. {\bf 127} 
(1999), no. 1, 33--53.

\bibitem[DP]{dec-pro}
C. De Concini, C. Procesi, 
Nested sets and Jeffrey-Kirwan residues. 
Geometric methods in algebra and number theory, 139--149, 
Progr. Math., 235, Birkh\"auser Boston, Boston, MA, 2005. 

\bibitem[De]{den-os}
G. Denham, The Orlik-Solomon complex and Milnor fibre 
homology. Topology Appl. {\bf 118}  (2002),  no. 1-2, 45--63. 

\bibitem[DP]{dim-pap}
A. Dimca, S. Papadima, 
Hypersurface complements, Milnor fibers and higher homotopy groups of 
arrangements. 
Ann. of Math. (2) {\bf 158}(2003), no. 2, 473--507. 

\bibitem[ESV]{esv}
H. Esnault, V. Schechtman, E. Viehweg, 
Cohomology of local systems on the complement of hyperplanes.
Invent. Math. {\bf 109} (1992), no. 3, 557--561. 

\bibitem[JO]{jew-orl}
K. Jewell, P. Orlik, 
Geometric relationship between cohomology of 
the complement of real and complexified arrangements. 
Top. and its appl. {\bf 118} (2002), 113--129. 

\bibitem[K]{koh-hom}
T. Kohno, 
Homology of a local system on the complement of hyperplanes. 
Proc. Japan Acad. {\bf 62} Ser. A (1986), 144--147. 

\bibitem[OS]{orl-sol}
P. Orlik, L. Solomon, 
Combinatorics and topology of complements of hyperplanes. 
Invent. Math. {\bf 56} (1980), 167--189. 

\bibitem[OT]{orl-ter}
P. Orlik, H. Terao, 
Arrangements of Hyperplanes. 
Grundlehren Math. Wiss. {\bf 300}, Springer-Verlag, New York, 1992. 

\bibitem[PS]{ps-h}
S. Papadima, A. Suciu, 
Higher homotopy groups of complements of complex hyperplane arrangements. 
Adv. Math. {\bf 165} (2002), no. 1, 71--100.

\bibitem[R]{ran-mor}
R. Randell, 
Morse theory, Milnor fibers and minimality of hyperplane arrangements. 
Proc. Amer. Math. Soc. {\bf 130} (2002), no. 9, 2737--2743. 

\bibitem[Sa]{sal-top}
M. Salvetti, 
Topology of the complement of real hyperplanes in $C\sp N$. 
Invent. Math. {\bf 88} (1987), no. 3, 603--618.

\bibitem[SS]{sal-sett}
M. Salvetti, S. Settepanella, 
Combinatorial Morse theory and minimality of hyperplane 
arrangements. preprint 2007. 

\bibitem[STV]{stv}
V. Schechtman, H. Terao, A. Varchenko,  
Cohomology of local systems and the Kac-Kazhdan condition 
for singular vectors. J. Pure Appl. Algebra {\bf 100} 
(1995), 93--102.

\bibitem[Su]{suc-tra}
A. Suciu, 
Translated tori in the characteristic varieties of 
complex hyperplane arrangements. 
Topology Appl.  {\bf 118}  (2002),  no. 1-2, 209--223. 

\bibitem[VG]{var-gel}
A. Varchenko, I. M. Gel'fand, 
Heaviside functions of a configuration of hyperplanes. 
Functional Anal. Appl. 21 (1987), no. 4, 255--270

\bibitem[W]{whi-ele}
G. W. Whitehead, Elements of homotopy theory. 
Graduate Texts in Mathematics, {\bf 61}. 
Springer-Verlag, New York-Berlin, 1978. 


\bibitem[Y]{yos-min}
M. Yoshinaga, 
Hyperplane arrangements and Lefschetz's hyperplane section theorem. 
To appear in Kodai Math. Journal. 

\bibitem[Z]{zas-face}
T. Zaslavsky, 
Facing up to arrangements: Face-count formulas for 
partitions of space by hyperplanes. Memoirs Amer. Math. Soc. {\bf 154} 1975. 


\end{thebibliography}
\end{document}